\documentclass[reqno,11pt]{article}
\usepackage{amsmath,amsfonts,amssymb,euscript,theorem}

\numberwithin{equation}{section}
\newtheorem{theorem}{Theorem}[section]

{\theorembodyfont{\rmfamily} 
}
\newtheorem{lemma}[theorem]{Lemma}

\newcommand{\al}{\alpha}
\newcommand{\e}{\varepsilon}
\newcommand{\de}{\delta}
\newcommand{\ga}{\gamma}

\newcommand{\la}{\lambda}

\newcommand{\om}{\omega}

\newcommand{\Om}{\Omega}

\newcommand{\si}{\sigma}
\renewcommand{\th}{\theta}
\newcommand{\R}{\mathbb{R}}
\newcommand{\C}{\mathbb{C}}

\newcommand{\pd}{\partial}
\newcommand{\donothing}[1]{{}}

\newenvironment{proof}{{\bf Proof.}}{\hfill\fbox{}\par\vspace{.2cm}}

\begin{document}

\title{ Solitary Wave Solutions for the Nonlinear Dirac Equations }
\author{{}\\[0mm]
 Meijiao Guan\\
\small{guanmj@math.ubc.ca } \\[1.5mm]
Department of Mathematics, University of British Columbia,
\\Vancouver, Canada, V6T 1Z2 }

\maketitle


\renewcommand{\S}{\mathbb{S}}

\begin{abstract} In this paper we prove the existence and local uniqueness  of stationary
states for the  nonlinear Dirac equation
\[
i \sum_{j=0}^{3} \ga^j \pd_j \psi - m\psi + F(\bar{\psi}\psi)\psi
=0
\]
where $ m >0$ and $ F(s) = |s|^{\theta}$ for $ 1\leq  \theta < 2.$
More precisely we show that there exists $\e_0 > 0$ such that for
$\omega \in(m - \e_0, m), $ there exists a solution $ \psi(t,x) =
e^{-i\omega t}\phi_{\omega}(x), x_0 = t, x = (x_1, x_2, x_3),$ and
the mapping from $ \omega $ to $ \phi_{\omega} $ is continuous. We
prove this result by relating the stationary solutions to  the
ground states of nonlinear  Schr\"{o}dinger equations.
\end{abstract}




\section{Introduction}
\label{sec:intro}

A class of nonlinear Dirac equations for elementary  spin-$\frac 12$
particles (such as electrons) is of the form

\begin{equation}\label{eq:main}
 i\sum_{j=0}^3 \ga^j \pd_j\psi - m\psi +
 F(\bar{\psi}\psi) \psi= 0.
\end{equation}
Here $F: \R \to \R $ models the nonlinear interaction.
 $\psi: \R^4 \to \C^4$ is a four-component wavefunction,  and $ m$ is a positive number.
 $ \pd_j = \pd/{\pd x_j},$
  and $\ga^j$ are the $ 4 \times 4$ Dirac
matrices:
\[
\ga^0 = \left(\begin{matrix} I_2 & 0 \\
                       0 & -I_2
         \end{matrix}
         \right),
\ga^k  =\left(\begin{matrix} 0& \si^k \\
                       - \si^k & 0
         \end{matrix}\right), \quad  k =1, 2, 3
\]
where  $ \si^k$ are Pauli matrices:
\[
\si^1 = \left(\begin{matrix} 0 & 1\\
                       1 & 0
         \end{matrix}\right), \quad
\si^2 = \left(\begin{matrix} 0 & -i\\
                       i & 0
         \end{matrix}\right), \quad
\si^3 = \left(\begin{matrix} 1 & 0 \\
                       0 & -1
         \end{matrix}\right).
\]
We define
 \[ \bar{\psi} = \ga^0 \psi, \quad \bar{\psi}\psi = (\ga^0 \psi, \psi)
  = \sum_{i=1}^2 (\psi_i, \psi_i) - \sum_{i=3}^4(\psi_i, \psi_i)\]
where $ (\cdot, \cdot )$ is the Hermitian inner product in $\C^1.$


Throughout this paper we are interested in the case
 \begin{equation}\label{eq:nonlinearit}
F(s) = |s|^{\theta}, \quad 0 <  \theta < \infty.
\end{equation} The
local and global existence  problems  for  nonlinearity as above
have been considered in \cite{EV, MNO}. For us, we seek standing
waves (or stationary states, or localized solutions of
\eqref{eq:main}) of the form
 \[
 \psi(x_0, x) = e^{-i \omega t} \phi (x)
\]
where $ x_0 = t, x = (x_1, x_2, x_3).$ It follows  that $ \phi :
\R^3 \to \C^4$ solves the equation
\begin{equation}\label{eq:s}
i \sum_{j=1}^3 \ga^j \pd_{j}\phi - m\phi + \omega \ga^0 \phi +
F(\bar\phi\phi) \phi =0.
\end{equation}
Different functions $ F$ have been used to model various types of
self couplings. Stationary states of the nonlinear Dirac field with
the scalar fourth order self coupling (corresponding to $ F(s) = s$
) were first considered by Soler \cite{S1} proposing them as a model
of extended fermions. Subsequently, existence of stationary states
under certain hypotheses on $ F$ was studied by Cazenave and Vazquez
\cite{CV}, Merle\cite{M} and Balabane \cite{BCDM}, where by shooting
method they established the existence of infinitely many localized
solutions for every $0 < \omega < m.$ Esteban and S\'{e}r\'{e} in
\cite{ES}, by a variational method, proved the existence of an
infinity of solutions in a more general case for nonlinearity
\[ F(\phi) = \frac 12(|\bar{\phi}\phi|^{\al_1} +
b|\bar\phi\ga^5\phi|^{\al_2})\phi, \quad \ga^5 =
\ga^0\ga^1\ga^2\ga^3
\]
for $ 0< \al_1, \al_2 < \frac 12.$ Vazquez \cite{V} prove the
existence of localized solutions  obtained as a Klein-Gordon limit
for the nonlinear Dirac equation $(F(s) = s).$  A summary of
different models with numerical and theoretical developments is
described by Ranada \cite{Ra}.

None of the approaches mentioned above yield a curve of solutions:
the continuity of $ \phi$ with respect to $\om$, and the uniqueness
of $ \phi$ was unknown. Our purpose is to give some positive answers
to these open problems. These issues are important to study the
stability of the standing waves, a question we will address in
future work.

Following \cite{S2}, we study  solutions which are separable in
spherical coordinates,
\[
\phi(x) = \left (\begin{matrix}  g(r)    &  \left(\begin{array}{ll}  1 \\ 0  \end{array}\right)  \\
                                i f(r)   &  \left(\begin{array}{ll}  \cos
                                \theta \\ \sin \theta e^{i\Phi}
                                \end{array}\right)
            \end{matrix}\right)
\]
where $ r = |x|, ( \theta, \Phi)$ are the angular parameters and $
f, g$ are  radial functions. Equation \eqref{eq:s} is then reduced
to a nonautonomous planar differential system in the $r$ variable
\begin{equation}\label{eq:sy}
\begin{split}
f'+ \frac 2r f & = ( |g^2 - f^2|^{\theta} - (m-\omega)) g  \\
 g'  & = ( |g^2 - f^2|^{\theta} - (m +\omega))f.
\end{split}
\end{equation}
 Ounaies in \cite{O}  studied the existence of solutions for equation \eqref{eq:s} using
a perturbation method. Let $ \varepsilon = m - \omega.$ By a
rescaling argument, \eqref{eq:sy} can be transformed  into a
perturbed system
\begin{equation}\label{eq:ssch}
\begin{split}
 & u' + \frac 2r u - |v|^{2\theta} v + v -
  (|v^2 - \varepsilon u^2|^{\theta} - |v|^{2\theta}) v  =0, \\
  & v' + 2m u - \varepsilon( 1 + |v^2 - \varepsilon u^2|^{\theta} )u =
 0
\end{split}
\end{equation}
If $\varepsilon =0,$ \eqref{eq:ssch} can be related to
 the  nonlinear Schr\"{o}dinger equation
\begin{equation}\label{eq:sch}
- \frac {\Delta v}{2m} + v - |v|^{2\theta}v = 0,   \quad \quad u = -
\frac {v'}{2m}.
\end{equation}
It is well known that for $ \theta \in (0, 2),$ the first  equation
in \eqref{eq:sch} admits a unique positive solution called the
ground state $ Q(x)$  which is smooth, decreases monotonically as a
functions of $ |x|$  and decays exponential at infinity(see
\cite{P}, \cite{SS} and references therein). Let $U_0=(Q, -
\frac{1}{2m}Q'),$ then
 we want to continue $U_0$ to yield a branch
of bound states with parameter $\varepsilon$ for \eqref{eq:ssch} by
contraction mapping theorem.

Ounaies carried out this analysis for $ 0 < \th < 1$ and he claimed
that the nonlinearities  in \eqref{eq:ssch} are continuously
differentiable. But with the restriction $ 0 < \th < 1$ we are
unable to verify it. The term $|v^2 - \varepsilon u^2|^{\theta}$ has
a cancelation cone when $ v = \pm \sqrt \e u.$ Along this cone, the
first derivative of $|v^2 - \varepsilon u^2|^{\theta}$ is unbounded
for $ 0 < \theta < 1.$ But Ounaies' argument may go through for
$\theta \geq 1,$ which gives us the motivation of the current
research. However we can not work in the natural Sobolev space
$H^1(\R^3, \R^2).$ Since $H^1(\R^3) \hookrightarrow L^6(\R^3),$
 we lose regularity.  To overcome these difficulties,
 we want to consider equation \eqref{eq:ssch} in the Sobolev space $W^{1,p}(\R^3, \R^2), p > 2$ and
  $\theta \geq 1.$

To state the main result, we introduce the following notations. For
any $ 1 \leq p \leq \infty, L^p_r =L^p_r(\R^3) $ denotes the
Lebesgue space for radial functions on $\R^3.\;   W^{1,p}_r=
W^{1,p}_r(\R^3) $ denotes the Sobolev space for radial functions on
$\R^3.$  Let $ X^p_r =W^{1,p}_r \times W^{1,p}_r, Y_r^p =L^p_r
\times L^p_r.$ Unless specified, the constant $C$ is generic and may
vary from line to line. In this paper, we assume that $m = \frac
12,$ since after a rescaling  $\psi(x)= (2m)^{\frac
1{2\theta}}\Psi(2m x),$ equation \eqref{eq:main} becomes
\[
i \sum_{j=0}^3 \ga^j \pd_j \Psi - \frac 12 \Psi +
F(\bar\Psi\Psi)\Psi = 0 .
\]
We  prove the following results:

\begin{theorem}
Let $ \varepsilon = m - \omega.$
 For $ 1 \leq \theta < 2$ there exists $ \varepsilon_0 = \varepsilon_0(\theta) > 0$ and a unique solution of
 \eqref{eq:sy}  $ (f,g)(\varepsilon) \in \mathcal{C}((0, \varepsilon_0), W^{1,4}_r(\R^3, \R^2))$
satisfying
\[
\begin{split}
f(r) & = \varepsilon^{\frac {\theta +1}{2\theta}}(-  Q'(\sqrt{\varepsilon} r) +
 e_2(\sqrt{\varepsilon}r))\\
 g(r) &  = \varepsilon^{\frac 1{2\theta}}(Q(\sqrt{\varepsilon} r) + e_1(\sqrt{\varepsilon} r))
\end{split}
\]
with
\[
\| e_j\|_{W^{1, 4}_r} \leq C \varepsilon \quad \mbox{for some} \quad
C(\th)  > 0, j = 1, 2.
\]
\end{theorem}
{\bf Remark}: The necessary condition $|\omega| \leq m $  must be
satisfied in order to guarantee the existence of localized states
for the nonlinear Dirac equation (see \cite{V}, \cite{MM}).\\

The solutions constructed in Theorem 1.1  have more regularity. In
fact, they are classical solutions and have exponential decay at
infinity.
\begin{theorem}\label{thm:reg}
There exists $C(\e)  >0, \si(\e) > 0 $ such that
\[
\quad |e_j(r)| + |\pd_r e_j(r)|
 \leq C e^{- \si r} \quad j=1,2.
\]
Moreover, the solutions $ (f, g)$ in Theorem 1.1  are  classical
solutions
\[
f, g \in \bigcap_{2\leq p < +\infty} W^{2,p}_r.
\]
\end{theorem}

\noindent {\bf Remark.} From the physical view point, the nonlinear
Dirac  equation with $F(s) = s$ (Soler model) is the most
interesting. In fact,  Theorem 1.1, Theorem 1.2 are both true for
the Soler model. In fact, from \eqref{eq:ssch} one can find out that
$ (v^2 - \e u^2) - v^2 = - \e u^2$   which is Lipschitz continuous.
An adaption of the proofs of the above theorems  will yield:

\begin{theorem}\label{cor:soler}
For the Soler model $ F(s) =s,$ there is a localized solution of
equation \eqref{eq:s} satisfying Theorem 1.1 and Theorem 1.2.
\end{theorem}

Next we proceed as follows. In section 2, we introduce several
preliminary lemmas. In section 3, we give the proof of Theorem 1.1,
Theorem 1.2.

\section{Preliminary lemmas}
\label{sec:lemmas}

We list several lemmas which will be used in Section 3.
\begin{lemma}\label{le:ineq1}
Let $ g: \R \to \R $ be defined by  $ g(t) = |t|^{2\theta} t,
\theta
>0, $ then
\[ \left|g(a+ \sigma) - g(a) - (2\theta +1) |a|^{2\theta} \sigma
\right| \leq (C_1 |a|^{2\theta -1}  + C_2
|\sigma|^{2\theta -1})|\sigma|^2
\]
 where $ C_1, C_2$ depends on $\theta$ and  $ C_1 = 0 $ if $ 0 < \theta \leq
\frac 12.$
\end{lemma}
\begin{proof}
We may assume that $ a > 0$ in our proof. It is trivial if $
\sigma =0.$ So we assume that $ \sigma \neq 0.$ If $ a < 2|\sigma|,$
then $ |a + \sigma | < 3 |\sigma| $ and
\[
\begin{split}
& \left|g(a+ \sigma) - g(a) - (2\th +1) a^{2\theta} \sigma
\right|\\
& \leq|g(a + \sigma)| + |g(a)| + (2\theta +1) |a^{2\theta} \sigma| \\
& <  C_1 |\sigma|^{2\theta +1}.
\end{split}
\]
If  $
a \geq 2|\sigma|,$ then
\[
 a + \sigma \geq 2|\sigma| + \sigma \geq
|\sigma| > 0,
\]
so that
\[ g(a+ \sigma)=
 (a + \sigma)^{2\theta +1}
.\]
 Taylor's theorem gives
\[
\begin{split}
g(a + \sigma) - g(a) -(2\theta +1)a^{2\theta}\sigma = \frac 12
g''(\xi)\sigma^2
\end{split}
\]
where $\xi$ is between $ a + \sigma $ and $ a.$ Since $ g''(\xi) =
2\theta(2\theta +1){\xi}^{2\theta -1},$ if $ 2\theta -1 < 0,$ then
\[
|g''(\xi)| \leq C |\sigma|^{2\theta -1}.
\]
If  $ 2\theta -1 > 0,$ we have
\[
|g''(\xi)| \leq C \max\{ (a + \sigma)^{2\theta -1}, a^{2\theta
-1}\}  \leq C( |a|^{2\theta-1}
 + |\sigma|^{2\theta -1}).
\]
Hence we prove the lemma.
\end{proof}

\begin{lemma}\label{le:ineq2}
For any $ a, b \in \R,  \theta >0,$ we have
\[
\left || a - b|^{\theta} - |a|^{\theta}\right| \leq C_1
|a|^{\theta -1} |b| + C_2 |b|^{\theta}
\]
where $ C_1, C_2$ depends on $\theta$ and $ C_1 = 0 $ if $ 0 < \theta \leq 1.$
\end{lemma}
\begin{proof}
The proof is basically similar to that of  the lemma
as  above. It is trivial if $ b = 0.$ So we may assume that $ b \neq
0$ and $ a > 0.$ If $ a < 2 |b|,$ then
\[
\left || a - b|^{\theta} - |a|^{\theta}\right| \leq C(
|a|^{\theta} + |b|^{\theta}) < C|b|^{\theta}.
\]
On the other hand, if $ a \geq 2 |b|,$ then $  | a - b| \geq a - |b|
\geq |b|.$ So  by using the mean value theorem
\[
\left || a - b|^{\theta} - |a|^{\theta}\right| = \th |t|^{\theta - 1} |b|
\]
where $ t$ is between $ a - b $ and $ a$. If $ \theta - 1 >0,$
then
\[
 |t|^{\theta -1} \leq C (|a|^{\theta -1} + |b|^{\theta -1}),
\]
hence
\[
\left|| a - b|^{\theta} - |a|^{\theta}\right| \leq C_1 |a|^{\theta
-1} |b| + |b|^{\theta}.
\]
If $ \theta - 1 < 0,$ then $|t|^{\theta -
1} \leq C |b|^{\theta -1},$ so that  we conclude
\[
\left || a - b|^{\theta} - |a|^{\theta}\right| \leq C
|b|^{\theta}.
\]
The proof is complete.
\end{proof}

\begin{lemma}\label{le:ineq3}
For  any  $a, b, c  \in \R, $  if $  1 \leq \theta < 2 ,$  then
\[
\left||a+ b +c |^{\theta } -  |a+b|^{\theta} - |a+c|^{\theta} +
|a|^{\theta}\right|  \leq C (|c|^{\theta -1} + |b|^{\theta
-1})|b|,\] where $ C$ depends on  $ \theta.$
\end{lemma}
{\bf Remark.} This inequality is symmetric about $b, c,$ so the right hand side
can be equivalently replaced by $C (|c|^{\theta -1} + |b|^{\theta
-1})|c|.$ Without loss of generality, we assume that $ |b| \geq |c|$ in the following.\\
\begin{proof}
For simplicity, let
\[
L =|a+ b +c |^{\theta } -  |a+b|^{\theta} - |a+c|^{\theta} +
|a|^{\theta}.
\]
It is trivial for $ \theta = 1,$  since if $ |a| \geq 5 |b|  $ then
 $ L  = 0. $
If $ |a| \leq 5 |b|,$
\[
|L| \leq C (|b| +  |c|).
\]
So next we  consider $ \theta > 1.$
If  $ |a| \leq 5 |b|,$
by triangle inequality and Lemma \ref{le:ineq2}, we have
\[
\begin{split}
| L | & \leq C( |a + c|^{\theta -1} + |a|^{\theta - 1} +
|b|^{\theta-1}) |b| \\
& \leq C (|c|^{\theta - 1} + |b|^{\theta -1}) |b|.
\end{split}
\]
If $ |a| \geq 5 |b|,$  by using Taylor's theorem
\[
 |L| =  C|( a + t_1b  + t_2  c)|^{\theta -2 } |bc|.
\]
where $ t_1, t_2\in (0, 1)$ and
\[
| (a + t_1b  + t_2  c)| \geq |a| - 2|b| - |c|
\geq  |c|.
\]
So if $ 1 < \theta <  2,$ we have
\[
|L|\leq C |c|^{\theta -1}|b|.
\]
The proof is complete.
\end{proof}

\begin{lemma}\label{le:pre}
Let $ 2\leq p \leq \infty,   f: \R^3 \to \R$  be radial and bounded.
Suppose $ f_r + \frac 2r f \in L^p_{loc}, \frac 2r f \in L^p_{loc}.$
 If $f_r + \frac 2r f \in L^p ,$ then $\frac fr \in L^p $
and
\[\|\frac fr\|_{L^{p}} \leq C \|\pd_r f + \frac {2}{r}
f\|_{L^{p}}\]
\end{lemma}
\begin{proof}
We begin with $ p= \infty.$ Using integration by parts
\begin{equation}\label{eq:self}
r^2 f(r) = \int_0^r (\pd_{\rho} f + \frac 2\rho f) \rho^2 d\rho.
\end{equation}
Hence
\[
|r^2 f| \leq \|f_r + \frac 2r f\|_{L^{\infty}}
\int_0^r s^2 ds = \frac{r^3}{3} \|f_r + \frac 2r f\|_{L^{\infty }}
\]
which gives
\[\|\frac fr\|_{L^{\infty}} \leq C \|\pd_r f + \frac {2}{r}
f\|_{L^{\infty}}.\]
Next let us consider $ p=2.$ Let $ 0 < r_1 < r_2  < \infty.$ Denote $ D= \{ x\in
\R^3,0 <  r_1 < |x| < r_2\}$ and

\[
I = 2\pi^2 \int_{r_1}^{r_2} ( f_r + \frac 2r f)\frac fr  r^2dr.
\]
By H\"{o}lder inequality,
\[
I \leq C \|\frac fr\|_{L^2(D)} \| f_r + \frac 2r f \|_{L^2(D)}.
\]
On the other hand, we have
\[
I = \frac 32 \| \frac fr\|_{L^2(D)}^2 + \pi^2( r_2 f^2(r_2) - r_1
f^2(r_1)).
\]
Since $  r_2 f^2(r_2) > 0$ we have
\[
 \| \frac fr\|_{L^2(D)}^2 \leq C\left (\|\frac fr\|_{L^2(D)} \|
(f_r + \frac 2r f)\|_{L^2(D)} +  r_1 f^2(r_1)\right).
\]
Let $ r_2 \to \infty, r_1 \to 0,$
 we  obtain
 \[\|\frac fr\|_{L^{2}} \leq   C \|\pd_r + \frac {2}{r}
f\|_{L^{2}}.\]
 The intermediate case $ 2 < p < \infty$ is a
direct result of interpolation .
\end{proof}

\section{Proof of the main theorems}
\label{sec:proof}

Similar to \cite{O}, we  use a rescaling argument to transform
\eqref{eq:sy} into a perturbed system. Let $ \varepsilon = m -
\omega$ (remember $m = \frac 12$). The first step is to introduce
the new variables
\[
f(r) = \varepsilon^{\frac {\theta +1}{2\theta}} u( \sqrt\varepsilon
r), \quad g(r) = \varepsilon^{\frac 1{2\theta}} v(\sqrt\varepsilon
r) \] where $ (f, g)$ are the solutions of \eqref{eq:sy}. Then $ (u,
v)$ solve
\begin{equation}\label{eq:rescale}
\begin{split}
 & u' + \frac 2r u - |v|^{2\theta} v + v -
  (|v^2 - \varepsilon u^2|^{\theta} - |v|^{2\theta}) v  =0, \\
  & v' +  u - \varepsilon( 1 + |v^2 - \varepsilon u^2|^{\theta} )u =
 0.
\end{split}
\end{equation}
Our goal is to solve \eqref{eq:rescale} near $ \varepsilon =0.$
If  $\varepsilon = 0,$ \eqref{eq:rescale} becomes
\begin{equation}\label{eq:rescale1}
\begin{split}
 & u' + \frac 2r u - |v|^{2\theta} v + v =0 \\
 & v' +  u =0.
\end{split}
\end{equation}
This yields the elliptic equation
\begin{equation}\label{eq:ground}
\begin{split}
  - {\Delta v}+ v = |v|^{2\theta} v  ,\quad
 u  = -  {v'}
\end{split}
\end{equation}
It is well known that for $ 0 < \theta < 2,$ there exists a unique
positive radial solution $Q(x)= Q(|x|)$ of the first equation in
\eqref{eq:ground} which is smooth and exponentially decaying. This
solution called a nonlinear ground state. Therefore $ U_0=( -Q', Q)$
is the unique solution to \eqref{eq:ground} under the condition that
$v$ is real and positive. We want to ensure that the ground state
solutions $U_0 $ can be continued to yield a branch of solutions of
\eqref{eq:rescale}.

 Let
\[
v(r) = Q(r)  + e_1(r),  \quad u(r) = -  Q'( r ) + e_2(r).
\]
Substitution into \eqref{eq:rescale}   gives rise to
\begin{equation}\label{eq:new}
\begin{split}
& e_2'(r) + \frac 2r e_2(r) + e_1  - ( 2\theta + 1) Q^{2\theta}e_1
= K_1 (\varepsilon, e_1, e_2)\\
& e_1'(r) + e_2(r) = K_2(\varepsilon, e_1, e_2)
\end{split}
\end{equation}
where
\[
\begin{split}
 K_1(\varepsilon,  e_1, e_2) & =|Q+e_1|^{2\theta}(Q+e_1) -
(2\theta +1) Q^{2\theta}e_1 - Q^{2\theta +1} \\
 & +(|v^2 - \varepsilon u^2|^{\theta}- v^{2\theta}) v \\
 K_2(\varepsilon, e_1, e_2) & =  \varepsilon ( 1 + |v^2 -
\varepsilon u^2|^{\theta}) u.
\end{split}
\]
Define  $ L $ the first order linear differential  operator  $ L :
{X^{p}_r}\to Y^p_r $ by
\[
L \left(\begin{array}{ll}
 e_1 \\
 e_2 \end{array}\right)= \left( \begin{matrix} 1
- (2\theta + 1) Q^{2 \theta} & \pd_r + \frac 2r \\
                           \pd_r    &   1
\end{matrix}\right)
\left(\begin{array}{ll} e_1 \\ e_2 \end{array}\right).
\]
Then we aim to solve the equation
\begin{equation}\label{eq:sy2}
Le = K(\varepsilon, e)
\end{equation}
where $ e = (e_1, e_2)^T, K(\varepsilon, e) =(K_1,
K_2)^T(\varepsilon, e).$ Let $ I = (0, \si), \si > 0.$  We say $ e
(\varepsilon) $ is a weak $X^p$-solution to  equation \eqref{eq:sy2}
if $ e$ satisfies
\begin{equation}\label{eq:sy3}
e = L^{-1}K(\varepsilon, e)
\end{equation}
for  a.e. $\varepsilon \in I.$ $L$ is indeed invertible as we learn from  the following lemma.

\begin{lemma} \label{le:iso}
Let  $0 < \theta< 2,$ the linear differential operator
\[
L = \left( \begin{matrix} 1
- (2\theta + 1) Q^{2 \theta} & \pd_r + \frac 2r \\
                           \pd_r    &   1
\end{matrix}\right)
\]
 is an isomorphism from $X^p_r $ onto $ Y^p_r$ for $  2 \leq p \leq \infty.$
 \end{lemma}
\begin{proof}
First we prove that $L$ is one to one.  Suppose that there exist
radial functions $ e_1, e_2\in W^{1,p}_r $ such that
\[
  L \left(\begin{array}{ll}
         e_1 \\ e_2
   \end{array}\right) =0.
 \]
Then
\begin{equation}
- \Delta_r e_1 + e_1 - ( 2\theta + 1) Q^{2\theta} e_1 = 0, \quad
e_2 = - e_1'.
\end{equation}
It is well known  (see, eg. \cite{TS}) that $e_1 =0$  is the unique
solution in $ H^1.$

Next we prove that $L$ is onto. Indeed $L$ is a sum of an
isomorphism and a relatively compact perturbation:
\[
L =
\left( \begin{matrix} 1  & \pd_r + \frac 2r \\
                           \pd_r    &   1
\end{matrix}\right) +
\left (\begin{matrix} - (2\theta + 1) Q(r)^{2 \theta} & 0 \\
                           0  & 0
\end{matrix} \right)= \tilde{L} + M.
\]
$M$ is relatively compact because of the exponentially decay of the
ground state at infinity.
 So we only need to prove that $\tilde{L}$ is an isomorphism from
$ X_r^p$ to $Y_r^p,$ i.e.
for any  $(\phi_1, \phi_2) \in L^{p}_r \times L^p_r,$ there exist
$ (e_1, e_2)\in W^{1, p}_r \times W^{1, p}_r $ such that
\[
\tilde{L} \left(\begin{array}{ll} e_1 \\ e_2
\end{array}\right) = \left(\begin{array}{ll} \phi_1 \\ \phi_2
\end{array}\right).
\]
It is equivalent to solve
\begin{equation} \label{eq:iso1}
\begin{split}
 e_1 + (\pd_r  + \frac 2r) e_2   & =
\phi_1  \\
\pd_r e_1 + e_2 &= \phi_2
\end{split}
\end{equation}
and show that   $ e_1, e_2 \in W^{1,p}_r.$  By eliminating  $e_2$ we
know that $e_1$  satisfies
\begin{equation}\label{eq:int}
 (- \Delta_r + 1) e_1 = \phi_1 -
 (\pd_r + \frac 2r)\phi_2.
\end{equation}
 Define  $
G(x) = (4\pi)^{-1}|x|^{-1}e^{-|x|}.$ \eqref{eq:int} has the solution
\[
\begin{split}
e_1  & = G(x)* \left(\phi_1 -
 (\pd_r + \frac 2r)\phi_2\right)\\
  &= G(x)* \phi_1 + \pd_r G(x) * \phi_2.
\end{split}
\]
Here we have used the property of convolution and the fact
$(\pd_r + \frac 2r)^*f(r) = - \pd_r f(r)$ in $ \R^3.$
By Young's inequality and $ G, \pd_r G \in L^1(\R^3),$ we have
\[
\|e_1\|_{L^p} \leq \|G\|_{L^1} \|\phi_1\|_{L^p} + \|\pd_r
G\|_{L^1} \|\phi_2\|_{L^p_r}
\]
which implies
\[ e_1 \in L^p_r.\]
 Similarly $e_2$ satisfies
\[
( - \Delta_r + 1 + \frac 2{r^2})e_2 = \phi_2 - \pd_r \phi_1.
\]
Let
$
H (x)=\frac {x_3}{|x|}G(x),
$
then
\[
\begin{split}
e_2 & =H(x) *\phi_2  - H(x)*(\pd_r \phi_1)\\
& = H(x) *\phi_2 + (\pd_r H + \frac 2r H) * \phi_1 \in L^p
\end{split}
\]
since $ H, (\pd_r + \frac 2r)H \in L^1(\R^3).$

To improve the
regularities of $e_1, e_2,$ we go back to \eqref{eq:iso1}.
Since
\[
\pd_r e_1 = \phi_1 - e_2 \in L^p_r,
\]
we have $ e_1 \in W_r^{1, p}.$ Regarding the regularity of $ e_2,$
we know that
\[
(\pd_r + \frac 2r)e_2 = \phi_1 - e_1 \in L^p_r.
\]
By Lemma \ref{le:pre}
\[
\begin{split}
\|\pd_r e_2 \|_{L^p} & \leq C (\|\frac {e_2}{r}\|_{L^p_r} +  \|(\pd_r +
\frac 2r)e_2\|_{L^p_r}) \\
& \leq C \|(\pd_r +
\frac 2r)e_2\|_{L^p_r} =C \| \phi_1 - e_1\|_{L^p_r}.
\end{split}
\]
Hence we have  $ e_2 \in W^{1,p}_r.$
\end{proof}

Now we are ready to construct solutions of \eqref{eq:sy3} by using
the contraction mapping theorem.

{\bf Proof of Theorem 1.1.} To prove  Theorem 1.1, we prove there
exists $\varepsilon_1 > 0$ such that
 for every $ 0 < \varepsilon <  \varepsilon_1,$
there is a  unique solution to equation \eqref{eq:sy3}
\[
e = L^{-1}K(\varepsilon, e)
\]
in a small ball in $X_r^4.$  First we must ensure  that $
K(\varepsilon, e)$ is well defined in $Y^p_r$ if $ e \in X_r^p.$
Recall  that
\[
\begin{split}
 K_1(\varepsilon, e_1, e_2) & =|Q+e_1|^{2\theta}(Q+e_1) -
(2\theta +1) Q^{2\theta}e_1 - Q^{2\theta +1} \\
 & +(|v^2 - \varepsilon u^2|^{\theta}- v^{2\theta}) v \\
 K_2(\varepsilon, e_1, e_2) & =  \varepsilon ( 1 + |v^2 -
\varepsilon u^2|^{\theta}) u.
\end{split}
\]
Let us consider $K_1,$ the estimate for $K_2$ is similar.
Since
\[
|K_1(\varepsilon, e)|\leq C_{\varepsilon, \theta}( |v|^{2\theta+1}
+ |u|^{2\theta+1})
\]
where $ C_{\varepsilon, \theta} $ is a real constant depending on
$\e, \th,$  it suffices to show that $ ( |v|^{2\theta+1} +
|u|^{2\theta +1 }) \in L^p.$ By Sobolev's embedding
\[
W^{1, p}(\R^3) \hookrightarrow L^q(\R^3)
\]
for any  $q$ if $ p > 3.$ We choose $p = 4$ in the following.
 The same argument is
available for $K_2.$ From Lemma \ref{le:iso}, we know that $ L^{-1}
K \in X_r^p.$

Fix $ \delta, $  to be chosen later. Consider the
set
\[
\Omega = \{e\in X^4_r; \|e\|_{X^4_r} \leq \delta\},
\]
and suppose $ e \in \Om.$  We know that
\[
\|L^{-1}K(\varepsilon, e)\|_{X^4_r}  \leq C ( \|K_1(\varepsilon, e)\|_{L^4_r} +
 \|K_2(\varepsilon, e)\|_{L^4_r}).
\]
Let $ K_1(\varepsilon, e)= K_1^n(\varepsilon, e) +
K_1^s(\varepsilon, e)$ where
\[
K_1^n(\varepsilon, e) =|Q+e_1|^{2\theta}(Q+e_1) - (2\theta +1)
Q^{2\theta}e_1 - Q^{2\theta +1}
\]
and
\[
K_1^s(\varepsilon, e) =\left(|(Q+e_1)^2 - \varepsilon(-Q'
+e_2)^2|^{\theta} - |(Q + e_1)|^{2\theta}\right)(Q+e_1).
\]
Thus
\[
\begin{split}
\| K_1\|_{L^4_r} \leq  \|K_1^s\|_{L^4_r} + \|K_1^n\|_{L^4_r}.
\end{split}
\]
For $\|K_1^s\|_{L^4_r},$ let $ a = v^2 = (Q + e_1)^2, b =
\varepsilon u^2 = \varepsilon ( -Q' +e_2)^2 $ in  Lemma
\ref{le:ineq2}, then
\[
\begin{split}
\| K_1^s\|_{L^4}  & \leq C_{\theta} \varepsilon \left\|| Q + e_1|^{2\theta -1}|-Q' +e_2|^2 +
|-Q' +e_2|^{2\theta}|Q +e_1|\right\|_{L^4_r}  \\
& \leq C_{\theta} \varepsilon ( \|Q\|^{2\theta +1}_{W_r^{1,4}} +
\|e\|^{2\theta +1}_{X^4_r})\\
&  \leq
 C_{\theta}\varepsilon ( \|Q\|^{2\theta +1}_{W_r^{1,4}} + \delta)\leq  \delta/4
\end{split}
\]
if  $ \de \leq 1 $ and $ \varepsilon$ is small enough such that
\[
C_{\theta}\varepsilon ( \|Q\|^{2\theta +1}_{W_r^{1,p}} +
 \delta)\leq \frac{\delta}{4}.
\]
For $\|K_1^n\|_{L^4},$ let $ a = Q(r), \sigma = e_1$ in Lemma
\ref{le:ineq1},  then
\[
\begin{split}
 \|K_1^n\|_{L^4_r} & \leq C_{\theta}( \| Q^{2\theta -1} e_1^2\|_{L^4_r} + \| e_1^{2\theta
+1}\|_{L^4_r}) \\
&  \leq C_{\theta}(\|e_1\|^2_{W^{1,4}_r} + \|e_1\|^{2\theta +1} _{W^{1,4}_r})\\
& \leq  C_{\theta} (\delta^2 + \delta^{2 \theta +1})\leq
2C_{\theta}\delta^2 \leq \delta/4
\end{split}
\]
if $ \delta \leq  \frac 1{8C_{\theta}}.$
A similar argument can be applied to $ K_2$ (with similar condition
on $\e, \de$) to obtain that
\[
\|K_2(\varepsilon,e)\|_{L^p_r} \leq \frac {\delta}{4}.
\]
Hence we obtain
\[ L^{-1} K(\varepsilon, e) \in \Omega.\]\\
Next we want to show that for any $ e, f \in \Omega,$ and $ \delta ,
\varepsilon $ as above,
\[
\|L^{-1}(K(\e,  e) - K(\e, f))\|_{Y^p_r} \leq \frac 34 \|e -
f\|_{X_r^p},
\]
i.e. $L^{-1} K$ is a contraction mapping. We have
\[
|K(\e,  e) - K(\e, f)|\leq |K_1^n( e) - K_1^n( f)|+ |K_1^s( e) -
K_1^s( f)| + |K_2(e) - K_2(f)|.
\]
We  compute the r.h.s. term by term. After rewriting $K_1^n( e)
 - K_1^n(f),$
\[
\begin{split}
|K_1^n(e) - K_1^n (f)| & \leq  \left||Q+e_1|^{2\theta}(Q+e_1) - |Q
+ f_1|^{2\theta} (Q + f_1) - (2\theta +1) |Q + f_1|^{2\theta}(e_1
-
f_1)\right| \\
& + (2 \theta +1) \left||Q+f_1|^{2\theta}(e_1 - f_1) - Q^{2\theta}
(e_1 - f_1)\right| = D_1^n + D_2^n.
\end{split}
\]
For $D_1^n,$ let $ a = Q + f_1 , \sigma =e_1 - f_1$  and by use of
Lemma \ref{le:ineq1}, then
\[
D_1^n \leq C( |Q+f_1|^{2\theta -1}+ |e_1 - f_1|^{2\theta-1}) |e_1
- f_1|^2.
\]
By Sobolev embedding and H\"{o}lder inequality, we have
\[
\begin{split}
\|D_1^n\|_{L^4_r} &  \leq C_{\theta}( \|e_1 - f_1\|^2_{W^{1,4}_r}
+\|e_1 - f_1\|^{2\theta +1}_{W^{1,4}_r}) \\
& \leq C_{\theta}(\delta + \delta^{2\theta})\|e_1 -
f_1\|_{W^{1,4}_r} \leq \frac 18\|e_1 - f_1\|_{W^{1,4}_r}
\end{split}
\]
if $ \de < \frac {1}{16 C_{\th}}.$  Using Lemma \ref{le:ineq2}, we
find
\[
|D_2^n| \leq C_{\theta} ( Q^{2\theta -2}+ |2Qf_1 + f_1^2|^{\theta
-1})
 |2Qf_1 + f_1^2||e_1 - f_1|.
\]
Hence
\[
\|D_2^n\|_{L^p_r} \leq C_{\theta}\delta \|e_1 - f_1\|_{W^{1,p}_r}
\leq \frac 18 \|e_1 - f_1\|_{W^{1,p}_r}.
\]
Then let us study $ K_1^s(e) - K_1^s(f):$
\[
\begin{split}
 K_1^s( e) - K_1^s( f) &  =
\left(| ( Q + e_1)^2 - \varepsilon ( - Q' +e_2)^2|^{\theta} - |(Q +
 e_1)^2|^{\theta}\right) (e_1 - f_1) \\
& +  (|( Q + e_1)^2 - \varepsilon (
- Q' +e_2)^2|^{\theta} - |(Q+e_1)^2|^{\theta} )(Q +f_1)\\
&
 -  (| ( Q + f_1)^2 - \varepsilon ( - Q' +f_2)^2|^{\theta} -
|(Q+f_1)^2|^{\theta})  (Q + f_1).
\end{split}
\]
Notice that  the first  line in the r.h.s. is easy to estimate since
\[
\begin{split}
& \left\| \left(| ( Q + e_1)^2 - \varepsilon ( - Q'
+e_2)^2|^{\theta} - |(Q + e_1)^2|^{\theta} \right) (e_1 - f_1)\right
\|_{L^p_r}  \\ &  \leq C_{\theta} \varepsilon \left\|(| Q +
e_1|^{2\theta -2}|-Q' +e_2|^2 + |-Q' +e_2|^{2\theta})(e_1 -
f_1)\right\|_{L^p_r} \\ & \leq \frac 18\|e_1 - f_1\|_{W^{1,p}_r}
\end{split}
\]
for $ \e$ sufficiently small.  For the second and the third line,
let us define
\[
\begin{split}
E(e,f)  & = (|( Q + e_1)^2 - \varepsilon (
- Q' +e_2)^2|^{\theta} - |(Q+e_1)^2|^{\theta} )(Q + f_1) \\
&
 -  (| ( Q + f_1)^2 - \varepsilon ( - Q' +f_2)^2|^{\theta} -
|(Q+f_1)^2|^{\theta})  (Q + f_1).
\end{split}
\]
We discuss the contractive property for  two different  situations $
\theta > 1 $ and $ \theta = 1$ separately. For $ \theta > 1,$
 we use Lemma
\ref{le:ineq3}. Set $ a =( Q + f_1)^2, b = ( Q + e_1)^2 - ( Q +
f_1)^2, c = -\varepsilon( Q + f_2)^2$ (notice that $b, c$ can be
taken sufficiently small),  and rewrite $E(e, f)$ to get
\begin{equation}
\begin{split}\label{eq:estimate}
| E(e, f)|  & \leq \left| |a+b+c|^{\theta} - |a +b|^{\theta} -
|a + c|^{\theta} + |a|^{\theta}\right|\sqrt {|a|}  \\
& +  \left| |( Q + e_1)^2 - \varepsilon ( - Q' +e_2)^2|^{\theta} -
|( Q + e_1)^2 -
\varepsilon ( - Q' +f_2)^2|^{\theta}\right | \sqrt {|a|}  \\
 & \leq C_{\theta}(|b|^{\theta -1} + |c|^{\theta -1})|b| \sqrt {|a|}\\
 & +  \left| |( Q + e_1)^2 - \varepsilon ( - Q' +e_2)^2|^{\theta} -
|( Q + e_1)^2 -
\varepsilon ( - Q' +f_2)^2|^{\theta}\right |\sqrt {|a|}.
\end{split}
\end{equation}
where for the last line, we applied  Lemma \ref{le:ineq2}.
 We obtain
\[
\|E(e, f)\|_{L^4_r} \leq C_{\theta}(\varepsilon^{\theta -1} +
\delta^{\theta -1}) \|e - f\|_{W^{1,4}_r} \leq \frac 18\|e -
f\|_{W^{1,4}_r}
\]
for $ \e, \de$ sufficiently small.  Hence we have for $ 1 < \theta <
2,$
\[
\|K_1(e) - K_1(f)\|_{L_r^{4}} \leq \frac 12 \|e -f\|_{W^{1,4}}.
\]
Next we prove that  $ E(e, f)$ is contractive  for $ \theta = 1 $ directly.
 Lemma \ref{le:ineq3} can  not be used
 since $ |b|^{\theta -1} = |c|^{\theta -1}
= 1.$ In \eqref{eq:estimate}, if $ |a|\geq \max\{5|b|, 5|c|\},$ then
\[
|a +b + c|- |a +b| - |a + c| + |a| = 0.
\]
Thus
\[
\|E(e,f)\|_{L^4_r} \leq C_{\theta}\varepsilon \|e - f\|_{W^{1,4}_r}
 \leq \frac 14 \|e - f\|_{W^{1,4}_r} .
\]
Hence we only need to consider $ E(e, f)$ if $ |a|$ is small, i.e. if $ |a|
\leq 5\max\{5|b|, 5|c|\},$
\[
|a +b + c|- |a +b| - |a + c| + |a| \leq C(|b| + |c|).
\]
Simply assume that $ |c| \leq |b|,$ we have
\[
\|E(e, f)\|_{L^p} \leq C_{\theta} \delta^{1/2} \|e - f\|_{W^{1,4}}
 \leq \frac 14\|e - f\|_{W^{1,4}} .
\]
Therefore if $ \theta= 1,$
\[
\|K_1(\varepsilon,e) - K_1(\varepsilon, f)\|_{L_r^{4}} \leq \frac 12
\|e -f\|_{W^{1,4}_r}.
\]
Similarly, we can prove that
\[
\|K_2(e) - K_2(f)\|_{Y^4_r} \leq \frac 14 \|e - f\|_{W^{1, p}_r}.
\]
Note we can satisfy all the condition above by choosing $ \de =
C_{\th} \e$ and taking $\e$ sufficiently small. Then  the
contraction mapping theorem implies  $ L^{-1}K $ has a unique fixed
point $ e(\varepsilon) \in \Omega$ which is a weak solution of
equation \eqref{eq:sy2}. The continuity w.r.t. $\e$ follows from the
continuity w.r.t $\e$ of the map $L^{-1}K$ and its contractibility.
This completes the proof of Theorem 1.1. $\square$\\

Let us see why a solution of equation \eqref{eq:sy2} which is in $
X_r^{1, 4}$ has more regularity. This is done by using a standard
bootstrap argument and the following standard lemma:
\begin{lemma}\label{le:reg1}
Let $ F: \C \to \C $ satisfy $ F(0) =0,$ and assume that there exists $ \al \geq 0$ such that
\[
|F(v) - F(u)| \leq C(|v|^{\al} + |u|^{\al})|v -u| \quad \mbox{for all} \quad u, v \in \C.
\]
Let
\[
\frac 1r = \frac{\al}{p} + \frac 1q , \quad 1 \leq p, q, r \leq \infty.
\]
It follows that if $ u \in L^{p}, \nabla u \in L^q,$ then $ \nabla F(u) \in L^r $ and
\[
\|\nabla F(u)\|_{L^r} \leq C \|u\|_{L^p}^{\al} \|\nabla u\|_{L^q}.
\]
\end{lemma}
{\bf Proof of Theorem 1.2.} First we can prove that
\[
e_1, e_2 \in \bigcap_{4 \leq p < \infty} W^{2,p}_r.
\]
Recall that $ e_1, e_2$ satisfy
\[
\left(\begin{array}{ll}  e_1 \\ e_2
\end{array}\right) = L^{-1}
\left(\begin{array}{ll} K_1 \\ K_2
\end{array}\right)
\]
and $ L$ is an isomorphism from $ X_p^r \to Y_p^r.$
We know that
\[
|K_1|\leq C_{\varepsilon, \theta}( |v|^{2\theta+1}
+ |u|^{2\theta+1})
\]
and
\[
|K_2| \leq \varepsilon (|u| + |v|^{2\theta+1}
+ |u|^{2\theta+1})
\]
Since $ e_1, e_2 \in W^{1,4}_r(\R^3) \hookrightarrow
L^{\infty}(\R^3),$ then $ K_1 \times K_2 \in L^{p}_r \times L^{p}_r$
for $ p \geq 4.$ By Lemma \ref{le:iso}, we have
\[
(e_1, e_2) \in  \bigcap_{4 \leq p < + \infty} W^{1, p}_r \times
W^{1, p}_r .
\]
Next from Lemma \ref{le:ineq1}, Lemma \ref{le:ineq2} and Lemma \ref{le:ineq3}
\[
\begin{split}
|K( e_1, e_2) - K( f_1, f_2)| & \leq C( |Q + Q'|^{2\theta}+  |e_1|^{2\theta } + |e_2|^{2\theta}
+ |f_1|^{2\theta } + |f_2|^{2\theta })\\
& \quad \quad (|e_1 - f_1| + |e_2 - f_2|).
\end{split}
\]
So by Lemma \ref{le:reg1},
\[
\nabla K_1 \times  \nabla K_2 \in  \bigcap_{4 \leq p < \infty}
L^{\frac {p}{2\theta}}_r
 \times L^{\frac {p}{2\theta}}_r.
\]
This gives that
\[
 (e_1, e_2)\in \bigcap_{4 \leq p < \infty} W^{2, p}_r\times W^{2,p}_r
\]
and
\[
\|e\|_{ W^{2,p}_r \times W^{2,p}_r }  \leq C\e.
\]
Going back to equation \eqref{eq:sy2}, we know that $e_1, e_2 \in
W^{3, p}_r \subset \mathcal{C}^2.$ So $(f,g)$ are classical
solutions.

 Moreover we show that $e_1, e_2 $ have
exponential decay at infinity. We know $ e_1, e_2$ are classical
solutions and $ |e_1|, |e_2| \leq C \e$ by Sobolev's embedding
theorem. Taking derivatives in \eqref{eq:rescale1} and after tedious
computations we find
\begin{equation}
\left\{\begin{array}{ll}
e_1{''} - e_1 = \de_1(r) e_1 + \de_2(r) e_1' + \de_3(r) Q& \mbox {for $ r$  large } \\
e_2{''} - e_2 = \si_1(r)  e_2+ \si_2(r) e_2' + \si_3(r)Q & \mbox
{for $r$ large}
\end{array} \right.
\end{equation}
where $ \si_i, \de_i \in W^{2, p} $ and $ |\si_1|, |\de_1| \leq C \e
(i = 1, 2, 3)$ for $ r $ large.

We conclude  that there exist constants $ r_0, \nu(\e), C(\e)$
positive such that
\begin{equation}\label{eq:exp}
|e_1(r)| + |e_2(r)| \leq C e^{-\nu r} \mbox{ for  $r \geq r_0$}.
\end{equation}
We prove it by an application of the maximum principle. Without loss
of generality, suppose $e_1(r_0) = 2 \e $ ($ r_0$ is sufficiently
large). Let
\[ h(r) = e^{-\nu (r- r_0)} + \beta e^{ \nu (r - r_0)}
\]
where $\beta > 0 $ is arbitrary and $0< \nu <1  $ is to be
determined later. If $ g = e_1 - h,$ then $ g$ satisfies
\[
g'' = ( 1 + \de_1 ) g +  \de_2 g' + ( 1 - \nu^2 + \de_1) h
 + \de_2  h' + \de_3 Q
\]
Since $ h' = \nu(-e^{-\nu (r- r_0)} + \beta e^{ \nu (r - r_0)})\leq
\nu h$ and $Q \leq  h,$ then
\begin{equation} \label{eq:max}
g'' \geq ( 1 + \de_1)g + \de_2 g' + ( 1 - \nu^2 + \de_1 + |\de_3|) h
\end{equation}
with $ g(r_0) = e_1(r_0) - ( 1 + \beta) < 0, g(\infty) < 0.$ Thus we
claim that
\[
g(r) \leq 0 \quad  \mbox{for} \quad r \geq r_0,
\]
if $ \nu $ is small enough such that
\[
1 - \nu^2 + \de_1  + |\de_3| \geq  0.
\]
If the claim  is not true, then $g(r)$ obtains maximum at $ r= r_1$
and $ g(r_1) > 0.$ Thus $ g''(r_1) < 0, g'(r_1) = 0.$  But this
contradicts with equation \eqref{eq:max} since the right hand side
of \eqref{eq:max} is positive evaluated at $ r= r_1$.  Therefore the
claim is true if $ \nu \leq \sqrt{1 - C\e}$  and then
\[
e_1(r) \leq h(r) \quad \mbox{if $r$ is large enough.}
\]
Then similarly we can show that
\[
e_1(r) \geq - h(r) \quad \mbox{if $r$ is large enough.}
\]
Thus
\[
|e_1(r)| \leq h(r) = e^{-\nu (r- r_0)} + \beta e^{ \nu (r - r_0)}.
\]
Letting $ \beta \to 0$ , we have
\[
|e_1(r)| \leq C e^{-\nu r}.
\]
for $ r $ large enough. The  exponential decay estimate for $ e_2$
can be obtained in a similar way. Once we have \eqref{eq:exp}, it is
obvious that $ |\pd_r e_j(r)| \leq C e^{-\nu r}$ and $ e_j \in
{H^2}.$ This completes the proof of Theorem 1.2. $\square$
\\

\noindent {\bf Remark.} For $0 < \theta < 1,$ our method does not
work since Lemma \ref{le:ineq3} is not valid.  Let us consider a
special example. Suppose $e_2 = f_2 = 0,$ then
\[
E(e_1, f_1) = (|(Q+e_1)^2 - \varepsilon(Q')^2|^{\theta} -
|Q+e_1|^{2\theta} -|(Q+f_1)^2 - \varepsilon(Q')^2|^{\theta} +
|Q+f_1|^{2\theta})(Q +f_1)
\]
We want to know whether or not the following inequality is true
\begin{equation}\label{eq:small}
 |E(e_1(r), f_1(r))| \leq \frac 14|e_1(r) - f_1(r)|, \quad r
\in (0, \infty)
\end{equation}
 if $\varepsilon$ small enough. Letting $ r_0$ large enough and
  $  s  = \e^{\al}, \al > 0$ to be determined later,
 we assume that
\[
\begin{split}
 & Q(r_0) + e_1(r_0) = \sqrt\varepsilon |Q'(r_0)| ( 1 + s),\\
& Q(r_0) + f_1(r_0) = \sqrt \varepsilon |Q'(r_0)|.
\end{split}
\]
Then under this ansatz,
\[
\begin{split}
 |E(e_1(r_0), f_1(r_0))|& =[ (s^2 + 2s)^{\theta} - ((1+s)^{2\theta} -1)]
h^{2\theta+1}
 = g(s) h^{2\theta +1}, \\
 |e_1(r_0) - f_1(r_0)| & =   sh
\end{split}
\]
where $ h =\sqrt \varepsilon |Q'(r_0)| .$ Then
\[
|E(e_1(r_0), f_1(r_0))| = \frac{g(s)}{s} h^{2\th}|e_1(r_0) -
f_1(r_0)| .
\]
We claim that if $ \al > \frac {\th}{1-\th},$ then
\[
\frac{g(s)}{s} h^{2\th} \gg \frac 12, \quad \mbox{as} \quad \e  \to
0.
\]
In fact, we have
\[
g(s) \geq C s^{\theta}
\]
since
\[
(s^2 + 2s)^{\theta} \geq C s^{\theta}
\]
and
\[
| ( 1 +s)^{2\theta} -1| \leq C ( s + s^{2\theta}) \ll C s^{\theta}.
\]
So
\[
\frac{g(s)}{s} h^{2\th}  \geq C s^{\th -1} h^{2\th} = C
|Q'(r_0)|^{2\th} \e^{\th + \al(\th -1)}\gg \frac 12, \quad \mbox{as}
\quad  \e \to 0
\]
since $\th + \al(\th -1) < 0.$   The claim is proved and
consequently, \eqref{eq:small} does not hold
for every $ r \in (0, \infty).$\\

\noindent {\bf Acknowledgement.} The author should like to thank
professor Eric S\'{e}r\'{e} for bring this problem to our attention.
The author also would like to thank Tai-Peng Tsai and Stephen
Gustafson for their very useful discussions.



\small

\begin{thebibliography}{}




\bibitem{BCDM} \textsc{M. Balabane, T. Cazenave, A. Douady, F.
Merle}.  \textit{Existence of excited states for a nonlinear Dirac
field.}  Commun. Math. Phys. 119, 153-176 (1988)

\bibitem{CV} \textsc{T. Cazenave, L. V\'{a}zquez}.
\textit{Existence of localized solutions for a classical nonlinear
Dirac filed}. Commun. Math. Phys. 105, 35-47 (1986)


\bibitem{ES} \textsc{M. J. Esteban, E. S\'{e}r\'{e}}.
\textit{Stationary solutions of the nonlinear Dirac Equations: A
Variational Approach}. Commnu. Math. Phys. 171, 323-350 (1995)


\bibitem{EV} \textsc{M. Escobedo, L. Vega.}  \textit{ A semilinear
Dirac equation in $H^s(\R^3)$  for $ s> 1$.} SIAM J. Math. Anal. 2
(1997), 338-362,





\bibitem{M}\textsc{F. Merle}.  \textit{Existence of stationary
states for nonlinear Dirac equations}. J. Diff. Eq. 74(1), 50-68
(1988)





\bibitem{MM} \textsc{P. Mathieu, T. F. Morris.} \textit{Existence condition for spinor
solitions.} Phys. Rev. D V.30, No. 8, 1835-1836, 1984


\bibitem{MNO} \textsc{S. Machinhara, K. Nakanishi, T. Ozawa.}
\textit{Small global solutions and the nonrelativistic limit for the
nonlinear Dirac equation.} Rev. Mat. Iberoamericana  19 (2003), 179
-194


\bibitem{O} \textsc{H. Ounaies.} \textit{Perturbation method for a class of noninear Dirac equations}.
 Differential and Integral
Equations. Vol 13 (4-6), 707-720 (2000)

\bibitem{P} \textsc{S. I. Pohozaev}, \textit{Eigenfunctions of the
equation $ \Delta u + \la u =0$,} Soviet Math. Dokl., 5 (1965), pp.
1408¨C1411


\bibitem{Ra} \textsc{A. F. Ranada.} \textit{Classical
nonlinear Dirac field models of extended particles.} In: Quantum
theory, group, fields and particles (editor A. O. Barut). Amsterdam,
Reidel: 1982


\bibitem{S1} \textsc{M. Soler}. \textit{classical, stable nonlinear
spinor field with positive rest energy}. Phys. Rev. D1, 2766-2769
( 1970)


\bibitem{S2} \textsc{M. Soler}. \textit{classical electrodynamics for a nonlinear
 spinorfield: perturbative and exact approaches}. Phys. Rev.
 3424-3429 (1973)



\bibitem{SS} \textsc{C. Sulem, P.-L. Sulem.}
\textit{The Nonlinear Schr\"{o}odinger Equations: Self-Focusing and
Wave
 \; Collapse,} Springer-Verlag, Berlin, 1999.




\bibitem{TS}\textsc{S.-M. Chang, S. Gustafson, K. Nakanishi and
 T.-P. Tsai.}  \textit{
Spectra of linearized operators of NLS solitary waves.} SIAM Journal
on Mathematical Analysis 39 (2007), no 4. 1070--1111.


\bibitem{V} \textsc{L. Vazquez}.  \textit{Localized solutions of a nonlinear
spinor field}. J. Phys. A: Math. Gen., Vol. 10, No. 8, 1977(1361
-1368).





\end{thebibliography}
\end{document}